\documentclass[10pt]{article}
\usepackage{amsmath}
\usepackage{amssymb}
\usepackage{amsthm}
\usepackage{exscale}
\usepackage{csquotes}
\usepackage{epsfig}
\usepackage[mathscr]{eucal}
\usepackage{float}
\usepackage{bm}
\usepackage{bbm}

\newtheorem{theorem}{Theorem}[section]

\newtheorem{remark}[theorem]{Remark}

\newtheorem{definition}[theorem]{Definition}

\usepackage{graphicx}
\usepackage{subcaption}
\numberwithin{equation}{section}

\numberwithin{equation}{section}

\renewcommand{\phi}{\varphi}

\usepackage{xcolor}

\begin{document}
\title{Estimation by Stable Motions and its Applications}
\author{Anirban Das, Manfred Denker, Anna Levina and Lucia Tabacu}

\maketitle

{\it \small \qquad Dedicated to the memory of Wojbor\,A.\,Woyczy\'nski (1943--2021)}

\begin{abstract}
We propose a family of confidence intervals for nonparametric moment estimators if the  observations have large or 
infinite variances. The theoretical underpinnings which guarantee the soundness of the method are demonstrated. Extensive numerical simulations  show its superiority  over  bootstrap and normal approximation and its wide applicability.  Finally,  a confidence interval to estimate the coupling strength in neuronal networks is  proposed. 
\end{abstract}

\section{Introduction}\label{sec:1}

Resampling by using stable distributions was introduced in \cite{DDW}, however no application and no simulation about its performance was provided at that time. The purpose of this note is to fill this gap, to provide an effective algorithm for this resampling procedure and to apply the method for nonparametric  moment estimators. 

 Most existing methods for establishing asymptotically consistent confidence intervals of parameters rely on exact distributions or on normal approximation (the Central Limit Theorem (CLT)), provided the underlying statistics has finite variance. However, when the variance of the underlying distribution is infinite the CLT does not hold, no numerically efficient method seems to be known, again when the variance is large, the confidence intervals become  too wide to permit reliable conclusions. Thus the natural question here is to find new statistics not changing the parameter, but allowing for calculations  of asymptotic confidence intervals for moments. A class of such new statistics is  called resampling using stable motions in \cite{DDW}, here abbreviated as {\it stable resampling}. 
 
 The limiting distribution for our resampling method depends on stochastic integrals with respect to stable motions, introduced in \cite{RW} by Rosi\'nski and Woyczy\'nski, but  the quantiles of the limiting distribution are not directly calculable. It is possible to use almost sure versions for the convergence to the limiting distribution (abbreviated as  ASLT in analogy to the normal case) which are obtained in \cite{HKM}.
 
 We begin in Section \ref{sec:2} with a description of the program to be used in the simulations. It splits into two subroutines, the first is the estimation of confidence intervals, once quantiles are estimated. This is a standard formality, but relies on the limit theorem for stochastic integrals.
 The other novelty lies in the second subroutine which consists of  an estimation of the limiting distribution function using the ASLT method.
 In Section \ref{sec:3} the theoretical background is briefly sketched leading to the specific form used to calculate confidence intervals and to the formula which permits to estimate the unknown distribution from the data directly. This resembles vaguely to a bootstrap method, though it is almost surely and a different resampling method.
 
 Simulations are collected in Section \ref{sec:4}. It splits into several subsections. First, using a simple form of the ASLT-algorithm, simulations for the performance of the estimate of the unknown distribution function  are presented. 
 
 The second and major part of the section builds the core of  this note. Its applications are rather wide, due to the fact that it works for $U$-statistics in general.  In Section \ref{sec:4}, one specific class of $U$-statistic is used, to estimate confidence intervals for the mean of a distribution which has large or infinite variance. Its algorithm uses the full strength of the resampling algorithm. Clearly, other moment  estimators can be handled in the same way. The purpose of this subsection is to present simulation results
 of various types. First of all we  use distributions with heavy tails (power like distributions) since  we shall make use of such distributions in the last subsection. For these models (varying power law and perturbations thereof) we examine different statistics (parametrized by the order of the stable motion), different sample sizes and different confidence levels. In addition we compare the new method with the bootstrap method (established in \cite{DDW}) and approximation by normal distribution.
 
  Many processes observed in biology and neural science are subject to heavy tail distributions or those where estimation of the mean is affected by large variances even if sample sizes are large. The last subsection   provides such an example. We apply our algorithm, called the {\it Stable Resampling for Moments} (SRM)  to estimate the connection strength  in a complete neuronal network~\cite{EHE} from the expectation $E(X)$ of its avalanche size  distribution. It is shown that the method has advantages over the classical confidence interval estimation for asymptotically normal observables. This is due to the fact that the variance of the distribution is about $(E(X))^3$ and $E(X)\to\infty$ as the number of neurons increases.

\section{ Stable Resampling: The Algorithm}\label{sec:2}

This section contains the description of the algorithm on which our estimation procedure for moments is based upon. It will be called {\it Stable Resampling} and has a more formal  description which easily enables the transformation into a code: It consists of two subroutines, 
where the first one relies on the second one.\\

%\begin{algorithm}[H]\label{alg1}
\subsection{\bf Stable resampling for moments }\label{sec:2.1} %{alg_res}

\vspace{0.3cm}

This subroutine will be called {\it stable resampling for moments} and is abbreviated as $\text{SRM}(q,q',p,r_l,r_u,\delta)$. 

Let $q>p>1$, $q'<\frac qp$,  $r_u,r_l\in \mathbb N$, $10<n\in \mathbb N$ and $\delta\in \{0,1\}$. For each choice   of these parameters the SRM-algorithm to estimate the $q'$-th moment based on a sample of size $n$ proceeds as follows:\footnote{The requirement $10<n$ can be replaced by a similar requirement which prevents small number of observations having a big influence on the estimation of the limiting distribution.}

 \noindent{\bf Data:} {$X_1,X_2,\dots X_n$  iid sample following the distribution of $X$, with $E\bigg[\big|X\big|^{q}\bigg] < \infty$. }
 
\noindent{ \bf Result:} { Given $\alpha \in [0,1] $, the output is a one- or two-sided
$\alpha$-level confidence interval for $E\big[X^{q'}\big]$.}

 \noindent{\bf Sequential Steps};
\begin{enumerate}%[label{\roman*}]
    \item Calculate $\hat{\mu}$
    \begin{equation}
       \hat{\mu}\,:=\, \frac{1}{n} \sum_{i=1}^{n} X_i^{q'}
    \end{equation}
    \item Generate an  iid sample ${\bigg\{Y_i}\bigg\}_{i=1}^{n}$ (independent of the samples $X_1, \dots X_n$) from a stable distribution with location parameter $1$, skewing parameter $0$, stability parameter $p$, scale parameter $0.5$.\footnote{Again, the  scale parameter may be chosen differently for variants of the subroutine.}   
    \item For $i=1,...,n$ calculate 
    \begin{equation}\label{arg_asclt}
        W_i= X_i^{q'}Y_i - \hat{\mu}Y_i 
    \end{equation}
    Pass the sample ${\bigg\{W_i}\bigg\}_{i=1}^{n}$ as input to the $\text{ASLT}(p,r_l,r_u)$ in Section \ref{sec:2.2} to obtain two estimated  distribution functions $Edf_l$ and $Edf_u$. From $Edf_l$ calculate the lower $\frac{\alpha}{2}$ quantile $L$ and  from $Edf_u$ calculate the the upper $\frac{\alpha}{2}$ quantile $U$.
    \item $\overline{XY} \gets \frac{\sum_{i=1}^{n} {X_i}^{q'}Y_i}{n} $, $\, \, \overline{Y} \gets \frac{\sum_{i=1}^{n} Y_i}{n} $.
    \item If {$\delta == 1$}
{
    $C_U\gets \frac{\overline{XY}}{\overline{Y}} - \frac{n^{\frac{1}{p} -1} L}{\overline{Y}} $\;
    $C_L\gets \frac{\overline{XY}}{\overline{Y}} - \frac{n^{\frac{1}{p} -1} U}{\overline{Y}} $;
} \newline
 If {$\delta == 0$}
{
    $C_U\gets \frac{\hat{\mu}}{\overline{Y}} - \frac{n^{\frac{1}{p} -1} L}{\overline{Y}} $\;
    $C_L\gets \frac{\hat{\mu}}{\overline{Y}} - \frac{n^{\frac{1}{p} -1} U}{\overline{Y}} $\;
   
}
 \item Output the $\alpha$-level confidence interval $C(q,p,r_l,r_u,\delta)\, = \, (C_L,C_U)$ and the $\alpha/2$-level confidence intervals $(-\infty, C_U)$ and $(C_L,\infty)$.
\end{enumerate}
%\end{algorithm}

\vspace{0.3cm}

%\begin{algorithm}[H]
\subsection{\bf Estimating the p-re-sampled distribution function   }\label{sec:2.2}%{alg_ASCLT}

\vspace{0.2cm}

The two distribution functions needed for the SRN subroutine are called {\it p-resampled distribution  functions} and are calculated according to the following subroutine $\text{ASLT}(p,r_l,r_u)$.

\noindent{\bf  Data:} {$W_1,W_2,\dots W_n$ iid sample following the distribution of $W$, $n > 10$. 

\noindent{\bf  Result:} { A pair of distribution functions $Edf_u$ and $Edf_l$  }

 \noindent{\bf Initialization}\;
    $CL_{\text{vec}} \leftarrow$ Vector(size: $r_l$),
     $\, \, \, CU_{\text{vec}} \leftarrow$ Vector(size: $r_u$)\;
     
 \noindent{\bf Sequential Steps}\;
 $r_m \gets \max(r_l,r_u)$ \;\newline
 {\bf For} {$perm \gets1$ {\bf to} $r_m$ {\bf by} $1$} {\bf do}
    \begin{enumerate}%[label=(\roman*)]
        \item randomly permute  $W_1,W_2,\dots W_n$ to obtain $W'_1,W'_2,\dots W'_n$.
        \item Define $T_1, T_2,\dots T_{n-9}$ as
        \begin{equation}\label{eq:2.3}
            T_i= (i+9)^{-\frac{1}{p}}\sum_{i'=1}^{i+9} W'_{i'}
        \end{equation}
        \item For any $t$, define 
        \begin{equation}\label{eq:2.4}
            C(t)= \bigg(\sum_{i'=10}^{n}\frac{1}{i'} \bigg)^{-1}  \bigg(\sum_{i'=10}^{n}\frac{1}{i'} {\mathbbm{1}}_{(-\infty, t]} \big(T_{i'-9}\big) \bigg).
        \end{equation}
    {\bf If} {$perm \le r_l$}
    {
        $CL_{\text{vec}}\bigg[perm \bigg] \gets C$
    }\newline
   {\bf If} {$perm \le r_u$}
        {
    $CU_{\text{vec}}\bigg[perm \bigg] \gets C$
        }     
    \end{enumerate}

  \begin{equation*}
         Edf_l(t) = \frac{1}{r_l} {\sum}_{c \in CL_{\text{vec}}} c(t),  \; \; \;
         Edf_U(t) = \frac{1}{r_u} {\sum}_{c \in CU_{\text{vec}}} c(t).
  \end{equation*}

\section{Theoretical Justifications for the Stable Resampling for Moments}\label{sec:3}

We first demonstrate the theoretical underpinnings which guarantee the consistency of the method  $\text{SRM}(q,q',p,r_l,r_u,\delta)$ described in Subsection \ref{sec:2.1},  where $q'$ is the order of the moment being estimated. From \cite{DDW} (Theorem 3.3) we have:
\begin{theorem}\label{theo:3.1}
Given an iid sample $X_1,X_2,\dots, X_n, \cdots$  following the distribution of $X$,  any function $h$ satisfying $E\bigg[\big|h(X)\big|^{r}\bigg] < \infty$, and an iid sample $Y_1,Y_2,\dots, Y_n, \cdots$  following a centered $p$-stable distribution with $p$ satisfying $r>p>1$, and also $Y_i$ being independent of $X_{i'},\, \forall i,i'$,  then
\begin{equation}\label{eq:3.1}
   \frac{1}{n^{\frac{1}{p}}} \sum_{i \le n} \left(h(X_i)- E_{x \sim X}\big[h(x)\big] \right)Y_i \xrightarrow{\text{weakly}} G_{X,h,p},
\end{equation}
for some random variable $G_{X,h,p}$ whose distribution depends on $X, h$, and $p$.
\end{theorem}
The setting considered in Algorithm $\text{SRM}(q,q',p,r_l,r_u,\delta)$  is recovered if we apply the above theorem with $h(x)= x^{q'}$, such that $E\bigg[\big|X\big|^{q}\bigg] < \infty$ and  $p< \frac q{q'}$. \newline
Next we demonstrate that in Step (iii) of Algorithm \ref{sec:2.1}, the call to Algorithm \ref{sec:2.2} yields two distribution functions approximating $G_{X,h,p}$. This in conjunction with (\ref{eq:3.1}) will establish the veracity of  Algorithm \ref{sec:2.1}.

\begin{definition}\label{def:LED}
\textbf{p-resampled distribution:} For the setting considered in Algorithm $\text{SRM}(q,q',p,r_l,r_u,\delta)$, the p-resampled distribution  $G_{X,q',p}$ is defined as the unique limit (in the sense of converge in law) of $n^{1-1/p}\sum_{i=}^n W_i $ as $n$ tends to infinity, with $W_i$ defined in \eqref{arg_asclt}.
\end{definition}
From Theorem 4.1 of \cite{HKM} we get:
\begin{theorem}\label{theo:3.2}
Given $X_1,X_2,\dots X_n$ iid samples following the distribution of $X$, $h$ any real-valued function, $Y_1,Y_2,\dots Y_n$ iid samples following a $p$-stable distribution with mean $0$, and
$$
T_{n}(h) := \frac{1}{n^{\frac{1}{p}}} \sum_{i \le n} \left(h(X_i)- E_{x \sim X}\big[h(x)\big] \right)Y_i.
$$
If $Y_i$ is independent of $X_{i'},\, \forall i,i'$, and $T_{n}(h)  \xrightarrow{\text{weakly}} G $, for some random variable $G$, then for any $-\infty \, <a\, <b\,< \, \infty$ such that $a,b$ are continuity points of G, we have
\begin{equation}\label{eq:3.2}
    \lim_{n \to \infty} \frac{1}{\log n}\sum_{k=1}^{n} \frac{1}{k} {\mathbbm{1}}_{(a,b)}\bigg(T_{k}(h)\bigg) \, \to \, \mathcal{P}\bigg(G \in (a,b)\bigg)\, \, as.
\end{equation}
\end{theorem}

If Theorem \ref{theo:3.2}  holds true under the relaxed condition that $E[Y_i]=1$ (see proceeding paragraph), it implies that when $\text{SRM}(q,q',p,r_l,r_u,\delta)$  calls the subroutine $\text{ASLT}(p,r_l,r_u)$ in Step (c),  the two distributions which are returned both converge in distribution to the p-resampled distribution  $G_{X,q',p}$. 

 Theorem \ref{theo:3.2}  is  only applicable when $E[Y_1]=0$, we show it's veracity for non zero values of $E[Y_1]$ when $E_{x \sim X}\bigg[\big|h(x)\big|^{r}\bigg] < \infty$ for some $r>p>1$. Without loss of generality let $E_{x \sim X}\big[h(x)\big] =0$, then
\begin{align*}
T_{n}(h) &= \frac{1}{n^{\frac{1}{p}}} \sum_{i \le n} h(X_i)Y_i  \\
&=\frac{1}{n^{\frac{1}{p}}} \sum_{i \le n} h(X_i)(Y_i- E(Y_i))\, + \, \frac{E(Y_i)}{n^{\frac{1}{p}}} \sum_{i \le n} h(X_i)
\end{align*}
Theorem \ref{theo:3.2} is applicable for the first part of the sum above, the second part goes to zero because $E_{x \sim X}\bigg[\big|h(x)\big|^{r}\bigg] < \infty$, and $ r > p$ (see \cite{CS}).
\begin{remark}\label{rem:3.3}
  The sequence $\{T_i\}_{i=1}^{n-9}$ defined in (\ref{eq:2.3}) of Algorithm \ref{sec:2.2} are different for different permutations  $\{W'_i\}_{i=1}^{n}$. This enables different estimates of the quantiles to be derived by permuting the data. One can use this to reduce variance of the  final estimate by averaging across the estimates from different permutations. This is not possible for statistics which do not depend on the order of the samples. 
\end{remark}

\begin{remark} The results of \cite{DDW} and \cite{HKM} are also applicable to $U$-statistics in general. Some additional assumption has to be imposed. Since we are only simulating moment estimators we are not formulating this here.
\end{remark}
  
  \section{Simulation Results}\label{sec:4}
  
  As explained in the introduction, this section summarises our simulation results for the p-resampled distribution  in Subsection \ref{sec:4.1}, the  stable re-sampling for moments in Subsection \ref{sec:2.1} and an application to neuronal avalanches.
  
  \subsection{Simulations for estimating p-resampled distributions}\label{sec:4.1}

The Algorithm \ref{sec:2.2}  when called from Step (c) of Algorithm \ref{sec:2.1} returns estimates for the p-resampled distribution  (def \ref{def:LED}). Here we demonstrate the robustness of the distribution functions inferred by varying the number of samples, and comparing to estimates of the p-resampled distribution  obtained by bootstrapping. 
  
We took $N$ random samples $Z_i$  from a Pareto distribution with shape parameters $ 2$ and location parameter  $3$, and independently we took  $N$ samples  $Y_i$  from a stable distribution  with order $p=1.2$, shape $\gamma =1$, skewness $\beta=0$ and  mean $\delta=1$. We transformed the data of the first sample using the map 
$$ x\mapsto  f(x)=x\max\{\log |x|, 1), $$
in order to get a distribution not in the strict domain of attraction of a stable distribution (here called Pareto-like).
Then, for the sample
$$ X_i= f(Z_i),$$
define (in likeness of \eqref{arg_asclt} with $q'=1$)
\begin{equation}\label{sec4_1eq}
        W_i= X_iY_i - \hat{\mu}Y_i , \; \; \text{where}\; \; \hat{\mu}\,:=\, \frac{1}{n} \sum_{i=1}^{n} X_i. 
\end{equation}
The data was generated using the seed $1345$ in the R-software. We compute the empirical distributions produced by passing the sequence $W_i$ to $\text{ASLT}(1.2,1,1)$ (using using $T_1,T_2 \cdots,T_n$ instead of $T_1,T_2 \cdots,T_{n-9}$ in \eqref{eq:2.4}),  we also estimate the p-resampled distribution  by applying bootstrapping in two ways:
 1. bootstrapping the quantities $W_i= X_iY_i - \hat{\mu}Y_i ,$  where $ \hat{\mu}\,:=\, \frac{1}{1100} \sum_{i=1}^{1100} X_i$  and 2.  bootstrapping the quantities $W_i= X_iY_i - 12.565*Y_i $. Note $12.565$ is the  mean of the random variable $X_1$ calculated using $10^6$ samples.
 
 \begin{figure}
\centering
\begin{subfigure}[b]{.45\linewidth}
\includegraphics[width=\linewidth]{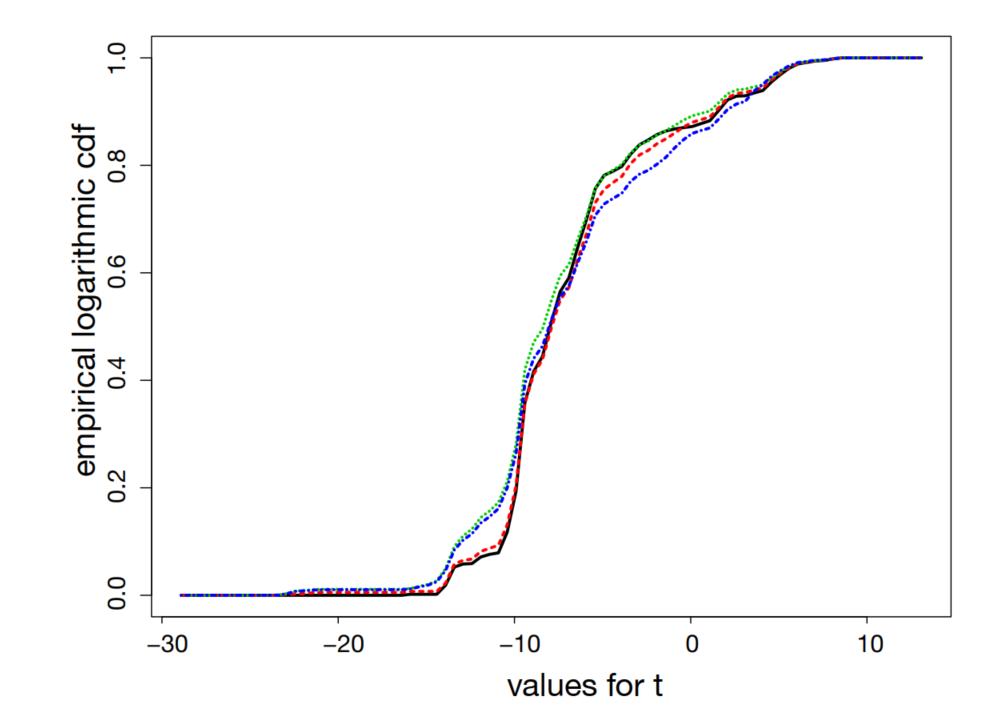}
\caption{ASLT}\label{Pareto}
\end{subfigure}

\begin{subfigure}[b]{.45\linewidth}
\includegraphics[width=\linewidth]{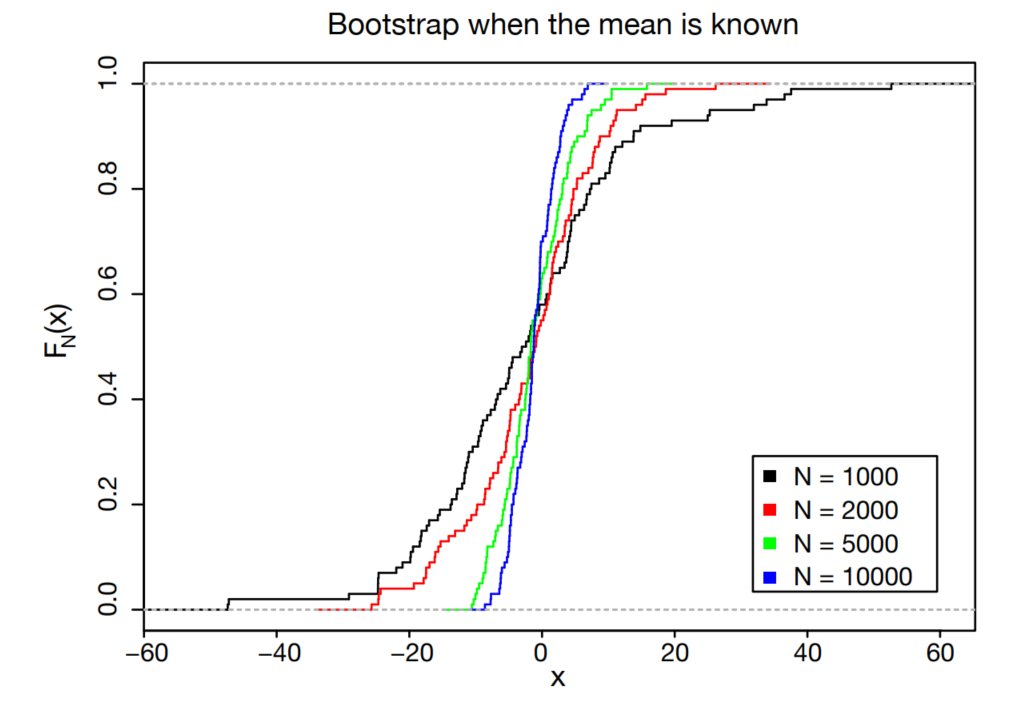}
\caption{Bootstrap (estimating mean)}\label{Paretoa}
\end{subfigure}
\begin{subfigure}[b]{.45\linewidth}
\includegraphics[width=\linewidth]{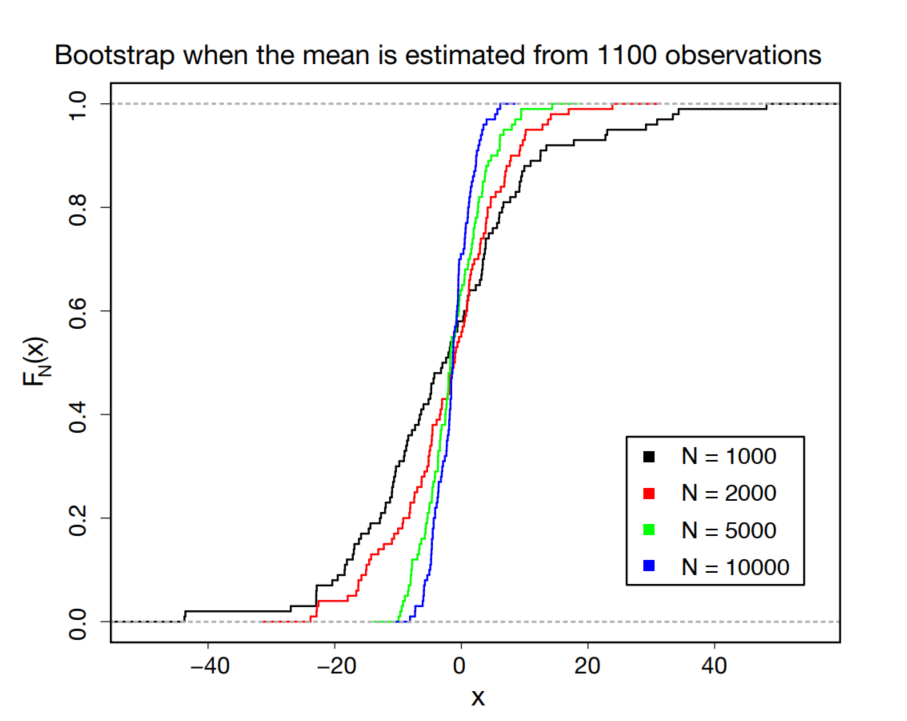}
\caption{Bootstrap (using the ground-truth value of the mean)}\label{Paretob}
\end{subfigure}
\caption{Estimating p-resampled distribution  for various sample sizes 1000(black), 2000(red), 5000(green) and
10000(blue), using different methods.}
\label{fig:animals}
\end{figure}

The distribution functions from applying Algorithm \ref{sec:2.2} clearly show the heavy tail behavior of distributions: Large  simulated values occur rarely, so are seen only in larger sample size simulations. The mean differs considerably from its median, hence the distribution is not symmetric around $0$.   The distribution does not have second moments but is in the domain of attraction of a normal. 
The graphics show that the estimation of the distribution function stabilizes quite well as $n\to \infty$ as expected from Theorem \ref{theo:3.2}.

Comparing the approximations in Fig \ref{Pareto} and Fig \ref{Paretoa}, it can be said that the ASLT approach is at least as good - if not better - than the bootstrap approximation. It also should be noticed that the bootstrap distribution seems to become  symmetric around $0$ and does not show the convergence pattern as the ASLT.  It is known that the ASLT approach can even be improved by using some permutations of the data and deleting some initial terms in the summation procedure.
This has been incorporated in the algorithms in Section \ref{sec:2}.

Figure  \ref{Paretob} shows the same approximation using bootstrap when the true mean $\mu$ is replaced by an estimated $\hat\mu$. The  graphics shows the same type of approximation, slightly shifted to the left, an effect due to the underestimation of $\mu$.

\subsection{Confidence intervals for the mean of power-like distributions}\label{sec:4.2}

Here we will investigate the performance of the SRM method for synthetic data-sets. First we introduce some metrics which will be used to evaluate the performance.

Recall that the cover probability
%  $\mathcal{M}$ be a method for generating an $\alpha$ confidence for a parameter $\theta$. Say a number of runs are used where the true value $\theta^{*} $ of $\theta$ is known, the cover probability
   is the proportion $p_c$ of runs where the unknown parameter $\theta^{*} $ lies in the $\alpha$-confidence interval divined by the estimation. If $p_c > \alpha$, then the method is called conservative. If  $p_c < \alpha$, then the method is termed permissive. 

Also the length of a confidence (CI length) is used to evaluate the accuracy of an estimation method.
 We take the average of the lengths of confidence intervals across all runs to obtain the CI length.\\

We will explicitly restrict ourselves to examples for estimating the first moment. Except otherwise stated, the data for  the examples studied in this section are generated from a  random variable $X \sim Z \log(Z)$, where $Z$ follows a Pareto distribution with scale parameter set to $1$, and shape parameter $q_{\text{P}}$ ($q_{\text{P}}$ is varied for various examples), these distributions are loosely called  \enquote{nearly power laws} . Given these settings $E[|X|^q]$ is finite for all $q < q_{\text{P}} $, and is $\infty$ for  all $q$ satisfying $q \ge q_{\text{P}} $. Whenever we apply $\text{SRM}(q,1,p,r_l,r_u,\delta)$, we are careful to choose a  value of $p$, satisfying $p < q_{\text{P}} $.

\vspace{0.8cm}

\noindent{\it \bf $\text{SMA}(q,1,1.2,r,r,\delta)$ versus normal approximation}\\

We compare the performance of $\text{SRM}(q,1,1.2,r,r,\delta)$, for both $\delta=0$, and $\delta=1$ to a standard CLT based method (which makes normal approximation) in a series of examples where we vary $ q<q_{\text{P}}$. We use the cover probability and length of confidence intervals to make comparisons. 500 runs of both methods were used, in each run 1000 samples were used, a two sided $95$ percentile interval was constructed, the exact results for the case $\delta=0$ are noted in Table \eqref{tab1}, the case $\delta=1$ showed similar results. For smaller values of $ q_{\text{P}}$ (underlying distribution is more heavily  tailed) the cover probability of normal approximations is very poor. At  $ q_{\text{P}}=4$, the performance of the normal method is excellent (Table \ref{tab1}). The cover probability produced by $\text{SRM}(q,1,1.2,r,r,0)$ is near the required .95 mark for all values of $ q_{\text{P}}$, similar trends hold for $\text{SRM}(q,1,1.2,r,r,1)$ .

\begin{table}[H]%[h!]
\caption{Cover probability ($P$) and CI length ($L$) of  $\text{SRM}(q,1,1.2,r,r,0)$ for the nearly power laws ,  $P_{nor}$ and $L_{nor}$ are the  cover probability and CI length of confidence interval for methods using normal approximations.}
\centering
\begin{tabular}{l||*{4}{c}r}
               & $P$ &  L &  $P_{nor}$ & $L_{nor}$  \cr
\hline\hline
$q_{P}=1.5$, r=5 & 0.948&  72.85  & 0.582 & 11.04 \\   
\hline
$q_{P}=1.6$, r=5 & 0.968 &38.06  & 0.648 &  4.72  \\
\hline
$q_{P}$=1.7, r=4  & 0.97 & 27.4 & 0.738 &  2.86\\  
 \hline  
$q_{P}$=2.1, r=2 & 0.946& 8.42 & 0.804 &  0.95 \\
 \hline
$q_{P}$=3.1, r=2  & 0.96 & 2.28 & 0.904 & 0.21 \\   
\hline 
$q_{P}$=4, r=2 & 0.97 & 1.33 & 0.942 &  0.11
\end{tabular}
\label{tab1}
\end{table}

\vspace{0.8cm}

\noindent{\it \bf Variation of the nearly power law}\\

We investigate the performance of the $\text{SRM}(q,1,1.2,r_l,r_u,0)$ method for different kinds of   \enquote{nearly power laws}. In particular we generate a Pareto distribution with scale parameter set to $1$, and shape parameter set to  $1.5$ ($q<1.5$), and pass it through a function $f$ to generate the ground truth data. The performance of the method is evaluated for various choices of $f$, see Table \ref{tab2}. The method gives acceptable performance for deriving both one sided and two sided confidence intervals for a variety of choices for $f$, which shows a type of robustness of the method with respect to perturbations.

\begin{table}[H]
\caption{$\mbox{P}(0)$ is the cover probability when $95\%$ symmetric two sided confidence interval is built for the mean, $\mbox{Pu}(0)$ is the cover probability when $97.5\%$ one sided upper confidence interval is built for the mean, $\mbox{Pl}(0)$ is the cover probability when $97.5\%$ one sided lower confidence interval is built for the mean, $n$ is the number of samples, $r_l$ and $r_u$  are the parameters used for the method $\text{SRM}(q,1,p,r_l,r_u,0)$, $q\sim 1.5$.  }
\centering
\begin{tabular}{l|l|l||c||*{4}{c}r}
               &function & $r_l$& $r_u$ & n & P(0) & Pl(0) & Pu(0) & L \cr
               \hline\hline
               p=1.2 & $x\log(|x|)$ & 5 & 5& 1000 & 0.958 & 0.978 & 0.98 & 43.82\\
               \hline\hline
               p=1.2 & $x/(1+\log(x))^{-2}$ & 2 & 2 & 1000 &0.93 & 0.992& 0.938 & 0.30\\
               \hline\hline
               p=1.2 & $x$ & 5 & 5 & 1000& 0.976 &0.978&0.998& 10.64\\
               p=1.2 & $x$ & 4 & 2 & 1000& 0.958  & 0.98& 0.978 & 9.06 \\
               \hline\hline
               p=1.2 & $x(1+\frac 12 \cos(x)$ & 5 & 5& 1000& 0.964 & 0.964 & 1 & 10.32 \\
               p=1.2 & $x(1+\frac 12 \cos(x))$ & 5 & 2 & 1000 & 0.968 &0.97 &0.998 & 12.93\\
               \hline\hline
               p=1.2 &$x\cos(x)$ & 5 & 5 & 1000 & 0.994 & 0.998 & 0.996 & 8.02\\
p=1.2 & $x\cos(x)$ & 2 & 2 & 1000 & 0.954 & 0.976 & 0.978 &  5.78  \\
\end{tabular}
\label{tab2}
\end{table} 

\vspace{0.8cm}

\noindent{\it\bf Using a single random permutation}\\

As noted in Remark \ref{rem:3.3} the ASLT part of the algorithm has the  ability to be able to reduce variance by taking permutations, here we demonstrate {\it taking only a single permutation} in Algorithm \ref{sec:2.1}  results in a rather poor final confidence interval(Table \ref{tab3}). 
\begin{table}[H]
\caption{Testing the performance of $\text{SRM}(1.5,1,p,r,r,\delta)$  for making $90 \%$ confidence intervals. $P(\delta)$ is the cover probability when $90\%$ symmetric two sided confidence interval is built for the mean, $Pu(\delta)$ ($Pl(\delta)$) is the cover probability when $95\%$ one sided upper (lower) confidence interval is built for the mean using $\text{SRM}(1.5,1,p,r,r,\delta)$, $n$ is the number of samples. The data is created by generating  a Pareto distribution with scale parameter set to $1$, and shape parameter set to  $1.5$, and passing the result  through a function $f(x) = x \log(x)$. 500 runs are made to derive the cover probabilities, for each run n samples are taken.  }
\begin{tabular}{l|l||c||*{7}{c}r}
               &$r$  & n & P(0) & Pl(0) & Pu(0)& P(1) & Pl(1) & Pu(1) &L \cr
\hline\hline
p=1.1 & 1& 1000& 0.804 & 0.894 & 0.91 & 0.816 & 0.958 & 0.858 & 38.92 \\
p=1.2 & 1& 1000& 0.74 & 0.872 & 0.868 &  0.76 & 0.956 & 0.804 & 24.7 \\
p=1.3 & 1& 1000& 0.702 & 0.882 & 0.82 & 0.698 & 0.957 & 0.74 & 35.03 \\
p=1.4 & 1& 1000& 0.676 & 0.874 & 0.802 & 0.662 & 0.964 & 0.698 & 14.27 \\
\end{tabular}
\label{tab3}
\end{table}

\vspace{0.8cm}

\noindent{\it\bf Performance under different stable distributions and sample sizes}\\
\begin{figure}[!]
 \includegraphics[scale=.4]{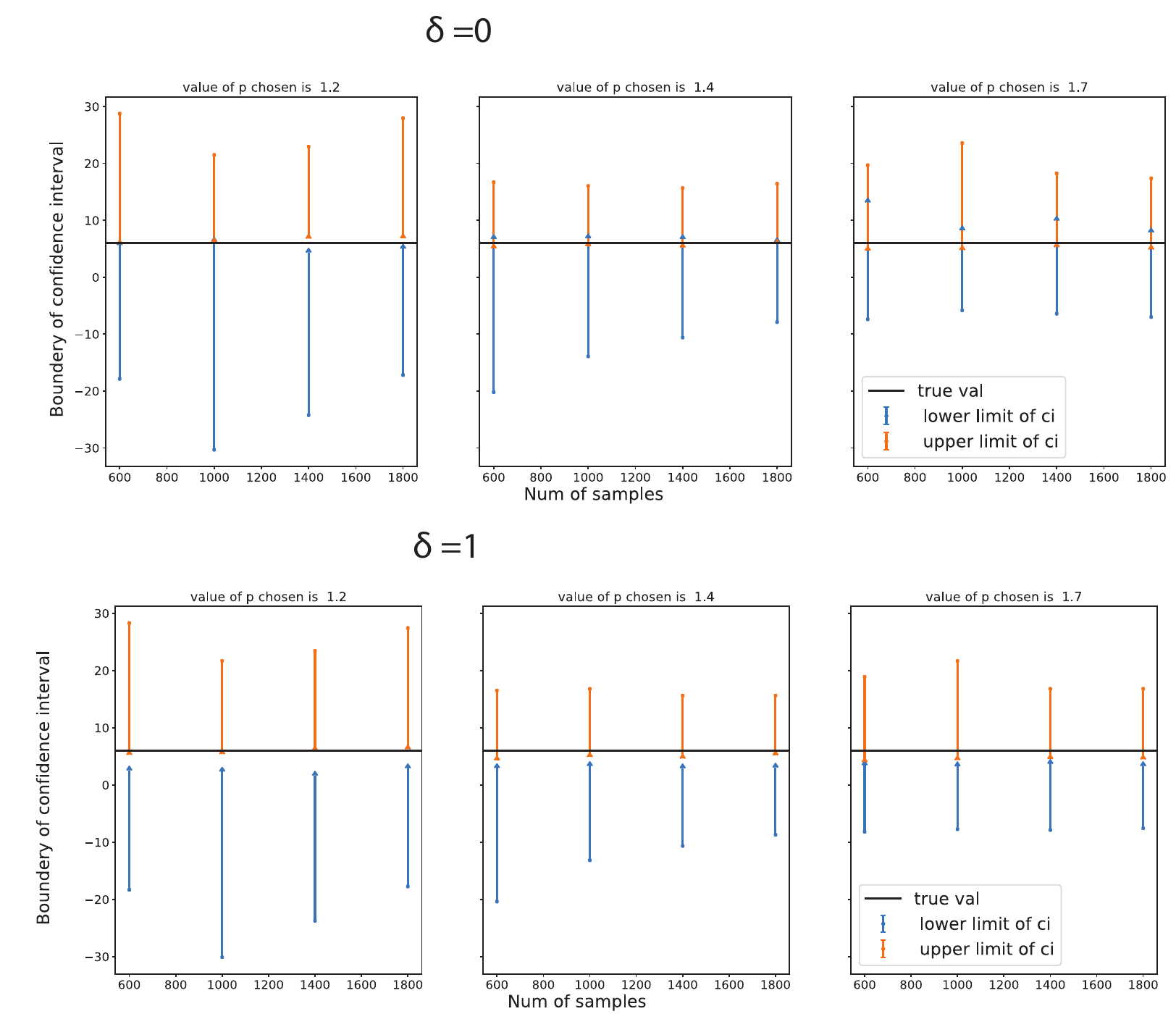}
 \caption{ 200 runs of $\text{SRM}(1,p,5,5,\delta)$ are made , with n samples used in every run. For each value of $n,\,p,\, \delta$, the mean (orange $\bullet$) and $2.5$-th percentile (orange $\bigtriangleup$) of the upper limit  of confidence intervals are noted. Also the the mean (blue $\bullet$) and $97.5$-th percentile (blue $\bigtriangleup$) of the lower limit  of confidence intervals are noted  }
\centering
\label{p_influence_onci}
\end{figure}
We demonstrate the performance of $\text{SRM}(q,1,p,5,5,\delta)$ in two simulations as p varies and $q\sim 1.5$. The data is generated by the same process as in Table \ref{tab3}. 200 runs are made (with $n$ samples used in every  run), for each run  $m$, the output from SRM is a confidence interval $(u(m), l(m))$. To visualize the performance, we observe the mean and $2.5$-th percentile of the upper limit  $u(m)$ (calculated across the  runs) in various controlled simulations where the value of $p$ and $n$ is varied. Similarly we also observe the mean and $97.5$-th percentile of the lower limit  $l(m)$ (see Figure \ref{p_influence_onci}). It is observed that for $p=1.2$, the $97.5$-th percentile of the lower limit and $2.5$-th percentile of the upper limit are close to the actual value of the parameter for all values of $n$. This clearly is not the case when $p=1.7$. The reason for the poor performance at $p=1.7$ is because the condition $p < q_{\text{P}} $ is violated for this example.

In the second simulation
we give the results for {SRM} under various choices of parameters when the underlying data is generated from the same nearly power laws as the data  for Table \ref{tab1}  (with $q_P=1.5$). The results are summarized in  Table \ref{tab4}.  It is notable that when $p$, $r_l$ and $r_u$ are fixed, the cover probability increases with the sample size indicating the method gets more conservative.

\begin{table}[H]
\caption{Performance of $\text{SRM}(q,1,p,r_l,r_u,\delta)$ when 500 runs of the method are used, each having $n$ samples. P($\delta$) denotes the cover probability of the two sided CI for the particular value of $\delta$, $\mbox{Pl}(\delta)$ ($\mbox{Pu}(\delta)$ ) is the cover probability of the one sided CI $[z_\alpha,\infty)$ (CI $(-\infty,z_\alpha]$) where $\alpha = 0.025$ ($\alpha = .975$).}

\centering
\begin{tabular}{l|l|l||c||*{7}{c}r}
               &$r_l$ & $r_u$ & n & P(0) & Pl(0) & Pu(0) & P(1)& Pl(1)& Pu(1)& L \cr
\hline\hline
p=1.1 & 5& 5&500 & 0.962 & 0.976   & 0.986 & 0.978  &  1& 0.978&74.77 \\   
p=1.1 & 5& 5&1000 &0.986  &  0.988 &0.998  & 0.99& 1&0.99& 89.69 \\ 
p=1.1 & 5& 5&2000 & 0.976& 0.98  & 0.996 & 0.996& 1& 0.996& 71.7 \\ 
\hline
p=1.1 & 2 & 5 & 500 & 0.914 & 0.932 & 0.982 & 0.952 & 0.99 & 0.962 & 74.34\\
\hline\hline
p=1.2  & 5& 5&500 & 0.944 & 0.98  & 0.964 &  0.948 &1& 0.948& 48.68\\   
p=1.2 & 5& 5&1000 & 0.958&  0.978 & 0.98  &  0.958& 1& 0.958& 43.82\\ 
p=1.2  &5 & 5&2000 & 0.974&  0.984 & 0.99 & 0.966&1&0.966& 48.67 \\ 
\hline
p=1.2 & 3 & 6& 500 & 0.932& 0.946 &0.964 & 0.946 &0.996 &0.95 &24.69 \\
p=1.2 & 6 & 7 & 500 & 0.934 & 0.97 & 0.968 & 0.97 & 1 & 0.97 & 50.45\\
p=1.2 & 3 & 6& 1000 & 0.93 &0.952 &0.978 & 0.96& 0.992 &0.968 & 52.52\\
p=1.2 & 6 & 7& 1000& 0.968 & 0.982 & 0.986  & 0.97 & 0.998 & 0.972 & 168.52\\
p=1.2 & 3 & 2& 2000 & 0.93& 0.972 &0.958 & 0.938&0.998 &0.94 & 236.21\\
p=1.2 & 2 & 4& 2000 & 0.946 &0.956 &0.99& 0.946&0.986 &0.96 & 41.4\\
\hline\hline
p=1.3  &5 &5 &500 & 0.91& 0.974  &0.936  & 0.888& 1& 0.888& 44.88 \\   
p=1.3 & 5& 5&1000 &0.936 & 0.964  & 0.972 & 0.928 & 0.998& 0.93& 37.03\\ 
p=1.3 &5 & 5&2000 & 0.954&0.972   & 0.982 &0.942 &0.998&0.944& 35.56\\ 
\hline
p=1.3 & 5 & 7 & 500 & 0.89 & 0.964 & 0.934 & 0.91 & 0.998 & 0.912 &  37.24 \\
 \hline\hline  
p=1.4  & 5& 5&500 & 0.858& 0.952  & 0.906 & 0.864 &0.998&0.866& 36.86\\   
p=1.4 & 5& 5&1000 & 0.898&  0.96 &0.938  &0.888 &0.994& 0.894& 43.51\\ 
p=1.4 & 5& 5&2000 & 0.926&0.964   & 0.962 &0.94 &0.996&0.944& 29.27 \\ 
\hline
p=1.4 & 6 & 7 & 500 & 0.88& 0.964 & 0.916 & 0.87 & 0.998& 0.972 &  31.86 \\
p=1.4 & 5 & 8 & 1000 & 0.928 & 0.96 & 0.968 & 0.92 & 0.998 & 0.922 &  34.02\\
 \end{tabular}
\label{tab4} 
\end{table}

\vspace{0.8cm}

\noindent{\it\bf Bootstrapping versus ASCLT-algorithm}\\

In the final set of simulations we see what happens if we pair the resampling approach with bootstrapping ~\cite{efr}, instead of the ASLT algorithm. Briefly, the bootstrap method $BRM(q',p,m)$ follows along the same lines of Algorithm \ref{sec:2.1}, with the only difference being that the two estimated  distribution functions $Edf_l$ and $Edf_u$ are not estimated using Algorithm \ref{sec:2.2} but by using a bootstrap sample of size $m$ (See \cite{Hall} for details). The results are given in table \ref{tab5}, the underlying data is generated from the same nearly power laws as the data  for Table \ref{tab1}  (with $q_P=1.5$). For all the methods studied two-sided symmetric  $95 \%$ confidence intervals for the mean are derived by making 500 runs, each using 1000 samples.   The cover probabilities for the $SRM$ methods remain relatively  closer to the desired $95\%$ for all values of parameters and sample size (see Table \ref{tab4} , section pertaining to $p=1.2$ for more details) when compared to the  BRM setting. Table \ref{tab5} shows that as m grows $\text{BRM}(1,1.2,m)$ becomes more conservative, similar trends were observed for $\text{SRM}(q,1,1.2,5,5,\delta)$ as the number of samples were increased, however the extent to which the results become conservative is far less for the $\text{SRM}(q,1,1.2, \cdot,\cdot,\cdot )$ methods. (see Table \ref{tab4}, section pertaining to $p=1.2$ for more details)

\begin{table}[H]
\caption{Comparison of different estimators for a 95\% CI for two-sided symmetric  $95 \%$ confidence intervals for the mean are derived by making 500 runs, each using 1000 samples}
\centering
\begin{tabular}{l|| c |c}
Estimator & Cover Probability & Length of CI \cr
\hline\hline
SRM(1.5,1,1.2,5,5,0) &               0.96  &  44.27  \\
RBM(1,1.2,50) &    0.984 &  33.72 \\
RBM(1,1.2,200) &  0.988 &  35.96  \\
RBM(1,1.2,500) & 0.994& 34.87 \\
Normal &  0.576&  6.26\\ 
\end{tabular}
\label{tab5}
\end{table}

\subsection{An application to neural avalanches}\label{sec:4.3}

The Abelian distribution is a distribution that is  important in models studying neural avalanches (see \cite{EHE} \cite{Levina} \cite{Levina2014abelian}), and it belongs to the class  of quasi Binomial II distributions (\cite{consul1975new}). Neural avalanches were observed in field studies, for example by Beggs et al. (\cite{Beggs2004}, \cite{Beggs2003}).  Cultured slices from the brain were attached to multielectrode ensembles and LFP (Local Field Potential) signals were recorded. The data retrieved showed  brief intervals of activity, when electrodes detected LFPs above the threshold. The period in between these  short bursts of activity was marked by idleness. A sequence  of such sustained
 activity was called an avalanche. There are models (\cite{EHE}) where the avalanche size (number of neurons firing during an avalanche) statistic follows an
 Abelian distribution. Recall that this
 distribution is a probability distribution on $\{1,2, \cdots ,N\}$ defined by the probability density
 	$$P(Z_{N,p} = b)= C_{N,p}{{N-1}\choose{b-1}}p^{b-1}(1-bp)^{(N-b-1)}b^{b-2},$$ where $C_{N,p}$ is the normalization constant defined by 
 	$C_{N,p}= \frac{1-Np}{1-(N-1)p} .$, where $N\in \mathbb N$ and  $p\in (0,\frac{1}{N})$ (\cite{Levina}, see also \cite{Levina2014abelian}). The $p$ in the Abelian distribution is often taken as $\frac{\alpha}{N}$, where $0<\alpha < 1$. It is known (\cite{Levina}, \cite{Levina2014abelian}}) that :
	\begin{equation}\label{abelian_mean}
	E(Z_{N,\frac{\alpha}{N}})= \frac{N}{N-(N-1)\alpha},\quad \text{hence}  \quad \lim_{N \to \infty} E(Z_{N,\frac{\alpha}{N}})= \frac{1}{1-\alpha}.
	\end{equation}
and (\cite{Das})
\begin{equation}\label{eq:ab_var}
\lim_{N \to +\infty} V(Z_{N, \frac{\alpha}{N}}) \,=\, \frac{\alpha}{(1- \alpha)^{3}}.
\end{equation}
We note that  (\ref{eq:ab_var}) can be proved borrowing results about Quasi Binomial II distributions (\cite{consul2006lagrangian}), and  asymptotic properties of incomplete gamma integrals. However, an elementary simple proof is given in \cite{Das}.

The parameter $N$ represents the number of neurons, in practice it is a large number. Also avalanches have been observed for collections of neurons of various sizes ~\cite{Beggs2003,Levina2017,Priesemann2013, tagliazucchi2012criticality, shriki2013neuronal,Petermann2009,das2019critical}, as such it is not presumed to be a phenomenon dependent on $N$. For a healthy brain the parameter $\alpha$ is hypothesized to be close to $1$ (\cite{EHE}). At $\alpha=1$ it is easy to show that (\cite{Levina}):
$$
\lim_{\alpha \to 1} \lim_{k \to \infty} \lim_{N \to \infty} \frac{P(Z_{N,\frac{\alpha}{N}} = k)}{C \;k^{-1.5}}= 1.
$$
All of this has two main consequences:
\begin{enumerate}
	\item For neural avalanche data the ratio of the underlying variance and mean will be very large.
	\item The distribution follows a nearly power law  with critical exponent $1.5$. This is in agreement with experimental observations where avalanche size distributions have been found to follow  power-law statistics, possibly with exponential cutoff~\cite{sornette2006critical, Beggs2003}.
	\item The quantity $\alpha$ is a useful parameter to be estimated from the data, since the extent of it's closeness to $1$ is thought to be a measure of the health of the brain. The quantity $\alpha$ can be estimated by estimating the mean.
\end{enumerate}
So this motivates us to estimate the mean of neural avalanche (using \eqref{abelian_mean}, one can estimate confidence interval for $\alpha$ using confidence intervals for the mean) data using the SRM algorithm. \\

\noindent{\bf Outline of simulations.}

\textbf{Data:} We use synthetic data. Our data is generated from a $1.5$ exponent  power law with upper cut-off at $x_m$ (we will analyze several data sets with different values of $x_m$). We will generate $n=1000$ iid instances of the data for each experiment, denote this by $X_1, X_2,\dots X_n$.\\

\noindent{\bf Results and discussion:}
\begin{figure}[htp]
	\centering
	\includegraphics[width=.3\textwidth, height=3.9cm]{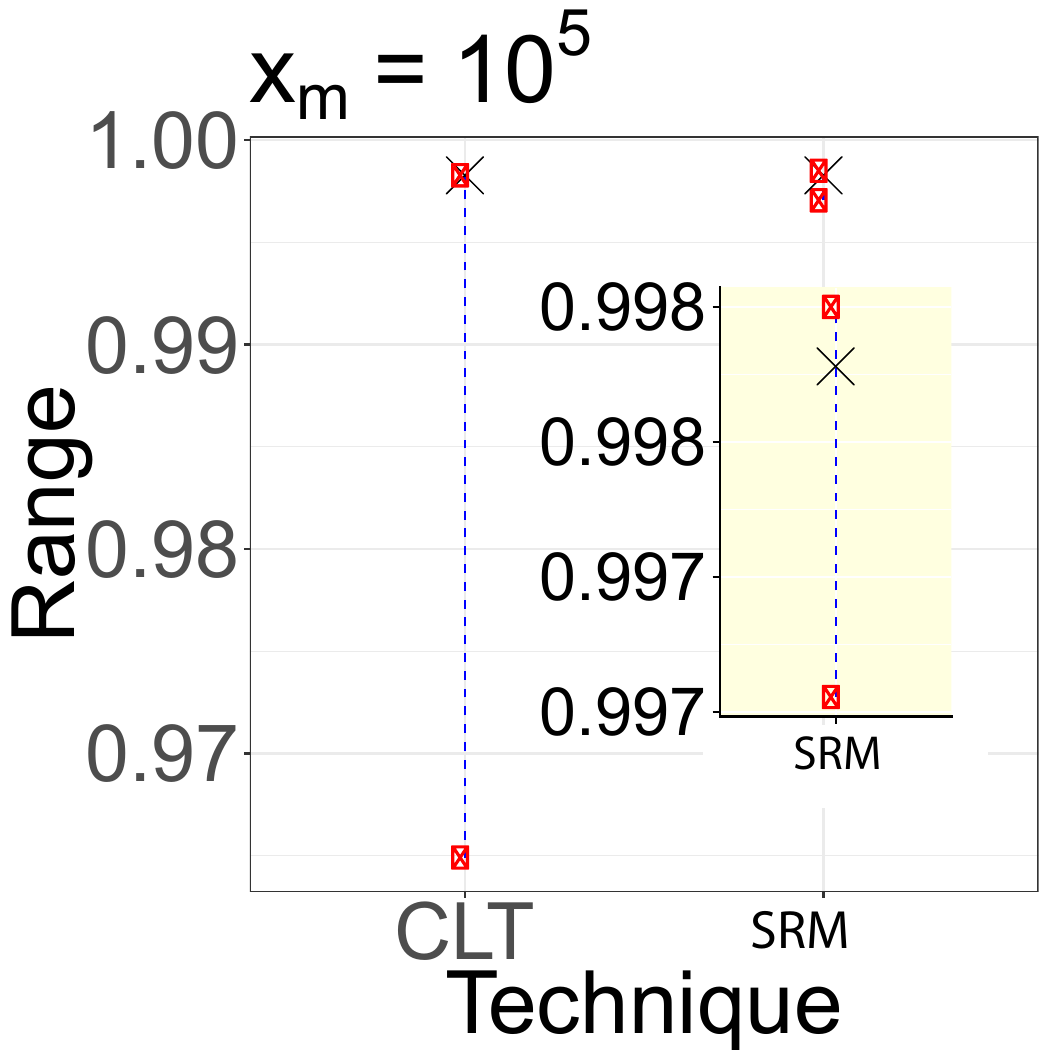}\hfill
	\includegraphics[width=.3\textwidth, height=3.9cm]{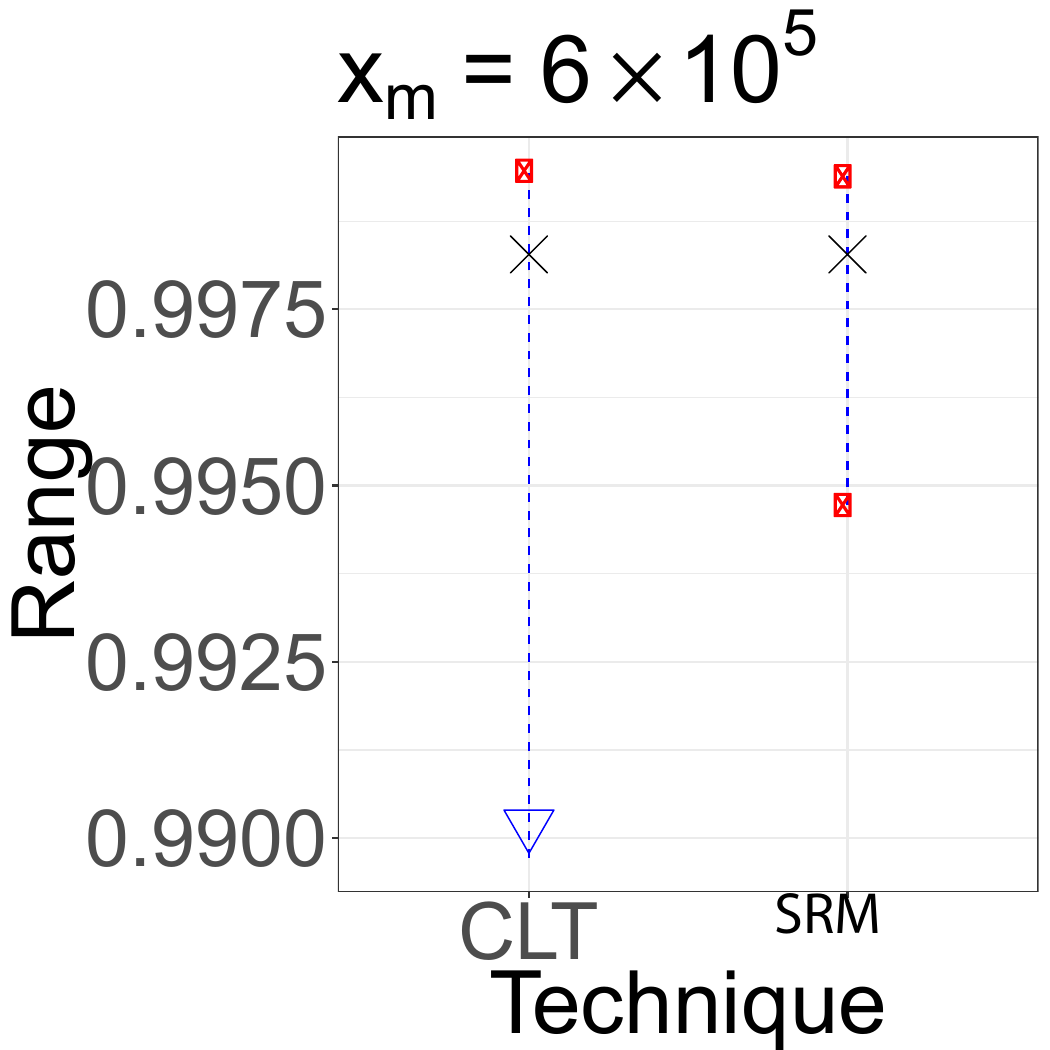}\hfill
	\includegraphics[width=.3\textwidth, height=3.9cm]{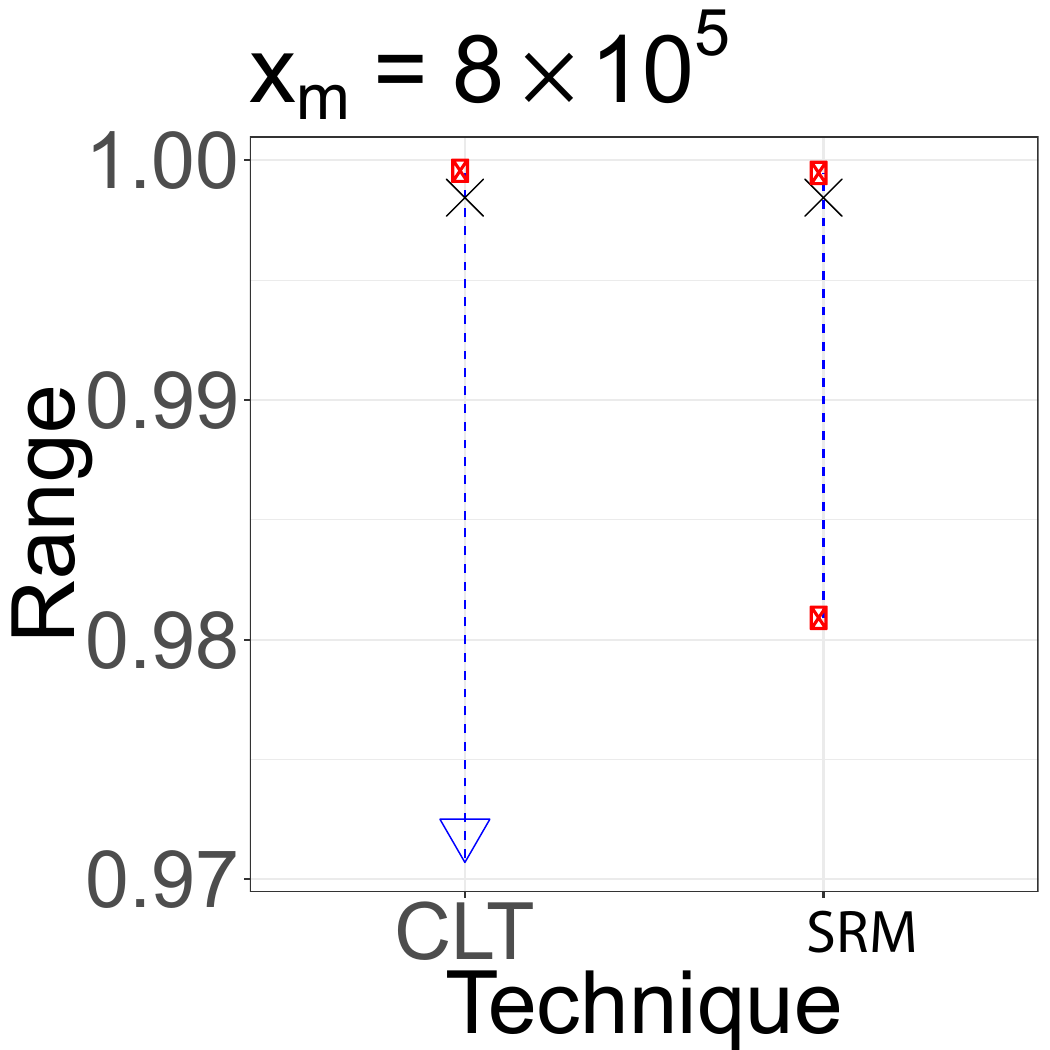}
	
	\caption{CLT and SRM method ($p$ value used is 1.7) for calculating confidence intervals for $\alpha$  for three different values of $x_m$: On the $x$-axes we indicate the method used to obtain confidence interval for $\alpha$. On the $y$-axes is shown the range of the  $4\%$ confidence interval obtained for each method. Red dots indicate the ends of the confidence intervals. The blue  $\bigtriangledown$ symbol indicates a lower bound for the confidence interval cannot be calculated using the method in question. To calculate confidence intervals we use $1000$ instances of synthetic data. The points indicated by $\times$ show the sample mean calculated from $900000$ instances of synthetic data. The inset on the leftmost  is to show the SRM results for this case more prominently.}
	\label{fig:htp}		
\end{figure}

Our simulation studies throw up some features worth noting:
\begin{enumerate}
	\item To check if our results are accurate we derive the sample mean from a much larger amount of synthetic data than what is used for establishing confidence intervals. This estimate will be called the precise sample mean and is marked by a $\times$  in  Figure \ref{fig:htp}.
	\item When $x_m=\infty$, the data is generated from a $1.5$ exponent power law over all of the positive integers. This distribution has both infinite first and second moments. In such a setting both the SRM and CLT methods will fail. As $x_m$ grows larger, the accuracy of both methods deteriorate. However because CLT method depends on higher moment conditions, it's accuracy deteriorates faster. Note, however for $x_m=  10^5$ the precise sample mean is near the center of the confidence interval calculated by the SRM method. But for $x_m=8 \times 10^5$ the lower bound of the confidence is quite far away from the precise sample mean.
	\item It is interesting to note that the methods can sometimes fail to give any lower bound for the confidence interval. The reason for this is as follows: We derive the confidence interval for $\alpha$ from the confidence interval for the mean $\mu$ using the understanding $\mu= \frac{1}{1-\alpha}$. For this one requires that both upper and lower confidence bound  estimates for the mean be positive, absence of  such conditions can result in lack of bounds. This happens in the case of the CLT method for $x_m=8 \times 10^5$ and $x_m=6 \times 10^5$. Although the underlying data is non-negative valued, the variance is so large that the lower confidence bound  obtained for the mean using CLT becomes negative.
\end{enumerate}

\vspace{0.5cm}
 
\noindent{\bf Acknowledgement:} This work is dedicated to the memory of Wojbor A. Woyczy\'nski, whose work laid the basis for the present investigation. The second named author was fortunate to discuss the finding with WAW in 2018/9 and we are thankful for his comments and advices.

Anirban Das\hfill{Manfred Denker}

Yale University \hfill{The Pennsylvania State University}

\vspace{1cm}

Anna Levina  \hfill{Lucia Tabacu}

University of  T\"ubingen \hfill{Old Dominion University}
\end{document}